\documentclass[12pt]{article}
\usepackage[utf8]{inputenc}
\usepackage{mathrsfs}
\usepackage{amsfonts,amssymb,amsmath,amsthm}
\usepackage{mathtools}
\usepackage{enumerate}
\usepackage{xcolor}
	\usepackage{amsfonts}
\usepackage{amsmath}
\usepackage{mathrsfs}

	\theoremstyle{definition}
	\theoremstyle{definition}

	\newcommand{\be} {\begin{equation}}
	\newcommand{\ee} {\end{equation}}
	\newcommand{\bn}{\begin{align}}
	\newcommand{\en}{\end{align}}
	\newcommand{\bea} {\begin{eqnarray}}
	\newcommand{\eea} {\end{eqnarray}}
	\newcommand{\Bea} {\begin{eqnarray*}}
		\newcommand{\Eea} {\end{eqnarray*}}
	\newcommand{\pa} {\partial}

	\newcommand{\ga} {\gamma}
	\newcommand{\Ga} {\Gamma}
	\newcommand{\Om} {\Omega}
	\newcommand{\om} {\omega}
	
	\newcommand{\la} {\lambda}
	\newcommand{\si} {\sigma}

	\newcommand{\vp} {\varphi}
	\newcommand{\var} {\varepsilon}

	\newcommand{\CC} {\mathbb{C}}
	\newcommand{\RR} {\mathbb{R}}
	\newcommand{\NN} {\mathbb{N}}
	\newcommand{\HH} {\mathbb{H}}
	
	\newcommand{\GG}{\mathcal{G}}
	
	\newtheorem{theorem}{Theorem}[section]
	\newtheorem{Pro}[theorem]{Proposition}
	
	\newtheorem{lemma}[theorem]{Lemma}
	\theoremstyle{remark}
	\theoremstyle{definition}
	\newtheorem{remark}[theorem]{Remark}
	\newtheorem{example}[theorem]{Example}
	\newtheorem{definition}[theorem]{Definition}
	
	\numberwithin{equation}{section}
	\def\R{\mathbb R}
	
	\def\C{\mathcal{C}}
	
	\def\qed{$\Box$}

	\def\sqr#1#2{{\vbox{\hrule height.#2pt
				\hbox{\vrule width.#2pt height#1pt \kern#1pt
					\vrule width.#2pt}
				\hrule height.#2pt}}}
	\def\square{\sqr74}
	\def\qed{{\unskip\nobreak\hfil\penalty50\hskip1em
			\hbox{}\nobreak\hfil\square \parfillskip=0pt
			\finalhyphendemerits=0 \par\goodbreak \vskip8mm}}

	\def\XXint#1#2#3{{\setbox0=\hbox{$#1{#2#3}{\int}$}
			\vcenter{\hbox{$#2#3$}}\kern-.5\wd0}}

	\title{Geodetically convex sets in the Heisenberg group $\HH^n$}
\author{Jyotshana V. Prajapat\\
\textit{\small Department of Mathematics}\\
\textit{\small University of Mumbai, Vidyanagari}\\
\textit{\small Mumbai 400 098, India}\\
\textit{\small jyotshana.prajapat@mathematics.mu.ac.in}\\
  Anoop Varghese\\
\textit{\small Department of Mathematics}\\
\textit{\small SIES College of  Arts, Science and Commerce, Sion}\\
\textit{\small Affiliated to University of Mumbai}\\
\textit{\small Mumbai 410 210, India}\\
\textit{\small  anoopv@sies.edu.in }}
\date{}
\begin{document}

	\maketitle

\begin{abstract}
We classify the geodetically convex sets and geodetically convex functions on the Heisenberg group $\HH^n$, $n\geq 1$.
\end{abstract}

	\section{Introduction}
A	Heisenberg group is the set of $2n+1$ tuples given by $\HH^n=\{(x,y,t) | x=(x_1,...,x_n),y=(y_1,...,y_n)\in \RR^n, t\in \RR \} $ with group operation $*$ defined as 
 \begin{equation} (x,y,t) * ( \tilde x, \tilde y, \tilde t )=(x+\tilde x,y+ \tilde y,t+ \tilde t+2\sum_{i=1}^{n}\tilde x_i y_i-x_i \tilde y_i). \nonumber\end{equation} 
 It is a Lie group with basis vector fields 
 \[ X_i = \pa_{x_i} + 2 y_i \pa_t, \quad  Y_i = \pa_{y_i} - 2 x_i \pa_t, 1\leq i\leq n, \quad T = \pa_t\]
 where $\{ X_i, Y_i: 1\leq i\leq n\}$ satisfy the H\''{o}rmander's condition 
 \[ [X_i, Y_j] = -4T \delta_{ij}, \quad 1\leq i, j \leq n.\] 
 Thus the Lie algebra generated by the vector fields $\{ X_i(p), Y_i(p): 1\leq i \leq n\}$  is the tangent space at the point $p$ and we denote the {\em horizontal space at $p$ \/} as ${\mathcal H}_p := \mbox{Span}\{ X_i(p), Y_i(p): 1\leq i \leq n\}$. In this paper, we will classify the geodetically convex subsets of $\HH^n$, extending the results proved in \cite{Monti} to higher dimension Heisenberg group. 
 
  Various notions of a convex domain and convex function in Heisenberg group have been defined, which are not necessarily equivalent. The first definitions appeared in \cite{Manfredi} where they discussed the notions of group convexity, horizontal convexity and viscosity convexity or convex in sense of viscosity.  For $\Om \subseteq 	\HH^n$ and $u: \Om \to \R$ an upper semicontinuous function, we say that $u$ is {\em  convex in the viscosity sense} or {\em v-convex}  in $\Om$ if 
\be 
(D^2_h u )^* \geq 0 \ee in the  viscosity sense, i.e., if $p \in \Om$ and $\vp \in C^2(\Om)$
with $\vp(p) = u(p)$ and $\vp(q) \geq u(q)$ for all $q$ in a neighbourhood of $p$ we have 	$(D^2_h \vp)^* (p)\geq 0$. 

For a point $p \in \HH^n$, let $[p * h^{-1}, p * h]$ denote the line segment which lies in the horizontal space $p * {\mathcal H}_o = {\mathcal H}_p$ which is obtained by left translation of the horizontal space at the origin, which will be denoted henceforth as $o $. An upper semicontinuous function $u$ defined on $\Om$ is said to be {\em horizontally convex } or just {\em $h$- convex} in $\Om$ if for all $p \in \Om$ and $h \in {\mathcal H}_0$ such that $[p * h^{-1}, p * h]  \in \Om$ we have \be u(p)  \leq \frac{ u(p * h) + u( p * h^{-1})}{2}. \ee
The above condition corresponds to {\em group convexity} if we require  $h \in \HH^n$ instead of  $h \in {\mathcal H}_o$. 
For more discussions, comparison and equivalences of various notions of convexity, we refer to \cite{Danielli}, \cite{Cristian}, \cite{Manfredi}, \cite{Capogna} .
	
Since straight lines are geodesics in Euclidean space, a more natural extension of convex  sets would be requiring that the function is convex along the geodesics joining any two points in its domain. This results in the definition of geodetically convex sets-\\
	\begin{definition}
		A subset $X$ of $\HH^n$ is said to be {\it geodetically convex} if every geodesic connecting every pair of points in $X$ lies in $X$.
	\end{definition}	
		
		It was proved in \cite{Monti} that the only geodetically convex sets in $\HH^1$  are empty set, points, arcs of geodesic and $\HH^1$. In this paper, we will extend this result to $\HH^n$, giving a proof which is true for all $n \geq 1$. Precisely,
		
\begin{theorem}\label{thm00} If $A$ is a geodetically convex subset of $\HH^n$ then either $A$ is an empty set, singleton set, an arc of a geodesic or the whole space $\HH^n$. 
\end{theorem}	
The proof of the Theorem \ref{thm00} is a consequence of  the following two theorems which will be proved in Section 	\ref{sets}.
	\begin{theorem}\label{main01} The smallest geodetically convex set containing two distinct points lying on a line parallel to the $t$ axis  is the whole space $\HH^n$.  
\end{theorem}	
\begin{theorem}\label{thm2} The smallest  geodetically convex subset of $\HH^n$ containing three distinct points not all of which lie on a geodesic is the whole space $\HH^n$.
\end{theorem}
Our proof relies on the description of the geodesics in all dimensions given by \cite{Zimmerman} and direct computations.  

We also classify here the geodetically convex functions on $\HH^n$.
\begin{definition}
	A function $u:\HH^n \to \RR$ is said to be {\it geodetically convex} if for every $p,q \in \HH^n$ and geodesic $\ga : [a,b] \to \HH^n$ from $p$ to $q$, the composite map $u \circ \ga : [a,b] \to \RR$ is convex(in the real sense).
\end{definition}
It was proved in \cite{Monti} that the only geodetically convex functions on $\HH^1$ are constant functions. We prove
\begin{theorem} \label{thm3.2}
	The only continuous geodetically convex function on $\HH^n$ is a constant function.
\end{theorem}

 Note that if we drop the continuity condition, we may have a geodetically convex function which is non constant.  The following example demonstrates that continuity is required by a geodetically convex function $u:\HH^n \to \RR$ for it to be constant. 
\begin{example}
	Consider the function $u:\HH^1 \to \RR$ defined as follows:
	\be u(x,y,t)= \left \{ \begin{array}{cc}
	1 & \text{if } t \notin L,\\
	0 & \text{if } t \in L
	\end{array} \right. \ee
	where $L=\{ (x,y,t) \in \HH^1: x=y, t=0 \}$.\\
Clearly, the function $u$ is not continuous in $\HH^1$ and also not constant. We claim that $u$ is geodetically convex. For if not,  there exists $ p_1, p_2 \in \HH^1$ and a geodesic $\ga: [0,1] \to \HH$ connecting $p_1$ to $p_2$ such that $u \circ \ga: [0,1] \to \RR$ is not convex (in Euclidean sense). This means, $\exists s_1, s_2 \in [0,1], \la \in (0,1)$ such that \be \label{3.2} u\circ \ga ( \la s_1 +(1-\la) s_2) > \la u\circ \ga ( s_1) +(1-\la) u\circ \ga(s_2).  \ee
This is possible only if \be \label{3.3} u\circ \ga ( \la s_1 +(1-\la) s_2) = 1 \ee
and \be \label{3.4} u\circ \ga ( s_1) = u\circ \ga(s_2) = 0. \ee
Using (\ref{3.4}) and definition of $u$, we have $\ga(s_1), \ga(s_2) \in L$. Therefore $\ga|_{[s_1, s_2]}=$ straight line segment from $\ga(s_1)$ to $\ga(s_2)$ lying in $L$. Since, $\ga(\la s_1 +(1-\la) s_2) \in \ga([s_1, s_2]) \subset L$, 
by definition of $u$,  $u\circ \ga ( \la s_1 +(1-\la) s_2)=0$ which contradicts $(\ref{3.3})$. Hence, the function $u$ is geodetically convex.
\end{example}

\begin{remark}
	In the above proof, we have used that the straight lines through the origin are geodesics and also the following result: 
	Let $p_1, p_2 \in  \RR^{2n} \times \{0\} \subset  \HH^n$ and $\ga: [a,b]\to \RR^{2n} \times \{0\} \subset \HH^n$ be a geodesic connecting $p_1$ to $p_2$. Let $\beta : [c,d] \to  \CC^n \times \{0\} \subset  \HH^n$ be a geodesic connecting $\ga(s_1)$ to $\ga(s_2)$ for some $s_1, s_2$ such that $a < s_1<s_2<b$. Then \[\ga([s_1, s_2]) = \beta([c,d]).\] 
\end{remark}

\begin{remark}
Possible generalizations of results proved in this paper to Carnot groups, groups of Heisenberg type and vertically rigid manifolds will be studied in a forthcoming paper. 
\end{remark}
Theorem \ref{thm3.2} is proved in Section \ref{functions}. We begin the next section by recalling prerequisites of Heisenberg  group required for our results.

\section{Prerequisites} \label{pre}

The Heisenberg group $\HH^n$ can also be written as $\CC^n \times \RR=\{(z,t) | z=(z_1,...,z_n)\in \CC^n, t\in \RR \}$ where the operation $*$ is defined as \be (z,t)*(\xi,s)= (z+\xi,t+s+2 \text{Im}\langle z, {\xi}\rangle) \nonumber\ee
where \[ \langle z, {\xi}\rangle =  \sum_{j=1}^{n}z_j\bar{\xi_j}\] is the Hermitian inner product, while it can be easily verified that \be \label{real}
(z, \xi) := \text{Im}\langle z, {\xi}\rangle = \sum\limits_{j=1}^n det [\xi_j, z_j] \ee is a skew symmetric bilinear form on $\R^{2n}$.
In the above expression, $\CC^{n}$ is identified wwith $\RR^{2n}$ using the map
\be (z_1, z_2, \ldots, z_n) =  (x_1+iy_1,...,x_n+iy_n) \mapsto    (x_1,...,x_n,y_1,...,y_n) . \nonumber \ee

Heisenberg group $\HH^n$ is a Lie group and its Lie Algebra($2n+1$ dimensional) of left invariant vector fields is generated by the following $2n+1$ vector fields:
\begin{align}
X_i&=\frac{\pa}{\pa x_i}+2y_i\frac{\pa}{\pa t} \;\text{for}\; i =1,2,...,n, \nonumber\\
Y_i&=\frac{\pa}{\pa y_i}-2x_i\frac{\pa}{\pa t} \;\text{for}\; i =1,2,...,n, \nonumber \\	
T&=\frac{\pa}{\pa t}.\nonumber \end{align} 
These vector fields satisfy the following commutation relations: 
\be [X_i,Y_j]=-4\delta_{ij}T, [X_i,T]=0, [Y_i,T]=0. \nonumber \ee
	
The {\em Horizontal Space }at a point $p$ is defined as  \[{\mathcal  H}_p  = \text{span}\{(X_1)_p,...,(X_n)_p,(Y_1)_p,...,(Y_n)_p \} \] and \[{\mathcal H}_p = \text{Ker}(\alpha_p) \;\; \forall p \in \HH^n \] where $\alpha$ is the standard contact form $ \alpha = dt + 2 \sum_{j=1}^n (x_j dy_j - y_j dx_j).$

 Instead of working with regular curves, we require that the regular curves that we consider have their tangent vector at every point lie in the horizontal space. Thus, an absolutely continuous curve $\ga : [0,S]\to \HH^n$ is said to be {\it  horizontal\/} if $\dot{\ga}(s)\in H_{\ga(s)}\HH^n$ for almost every $s \in [0,S]$. It is easy to verify (see \cite{Capogna}) that  an absolutely continuous curve  $\ga : [0,S]\to \HH^n$ given by $\ga = (x_1(s),...,x_n(s),y_1(s),...,y_n(s),t(s)) = (x(s), y(s), t(s))$ is horizontal iff \[ \dot{t}=2 \sum_{j=1}^{n} (\dot{x}_jy_j-\dot{y}_jx_j),\] i.e.,
 \be \label{e8.1} t(s)-t(0) = 2 \sum_{j=1}^{n} \int_{0}^{s} [\dot{x}_j(u)y_j(u)-\dot{y}_j(u)x_j(u)] du \ee for almost every $s\in [0,S]$. By Chow's theorem, 
	any two points in $\HH^n$ can be connected by a horizontal curve. Observe that a line segment  $\ga(s) =(1- s)p + s q $ joining two points $p= (x,y,0)$ and $q=(u,v,0)$ is horizontal iff
	\be \label{segment} ((x,y),(u,v))=
	\sum_{i=1}^{n}(u_i y_i-x_i v_i)=0. \ee
In case of $\HH^1$ this implies that the vectors $p$ and $q$ lie on the straight line segment passing through the origin. For $n \geq 2$, for a point $p \in \R^{2n}$,  a line $p +s \overrightarrow{w}$ in $\R^{2n}$ in the direction of vector $\overrightarrow{w}$ will be horizontal in $\HH^n$ iff 
	\[(\overrightarrow{w}, p) = 0\] i.e., for all $\overrightarrow{w}$ in the $2n-1$ dimensional  hyperplane $\{z\in \R^{2n}: (\overrightarrow{w}, p)   =0 \} \times \{t = 0\}$. Thus in higher dimensions, all these will further correspond to geodesics in the plane $t = 0$. In particular, any line segment joining  a point $p \in \R^{2n}$  to $-p$ will be horizontal and hence also a geodesic. 

Let $\Ga:[0,S]\to \HH^n$ be a horizontal curve and $\pi :  \HH^n \to \RR^n \times \RR^n \times \{0\}$ denote the projection map onto the first $2n$- coordinates, i.e., the $XY$ plane. Let $\ga = \pi(\Ga)$, the projection of $\Ga$ onto $\RR^n \times \RR^n \times \{0\}$. Then the $\HH^n$ length of $\Ga$ is defined as \be l_H(\Ga) = l_E(\ga)\ee where $l_E$ is the Euclidean length given by $l_E(\ga)= \int_0^S |\dot{\ga}|^2 ds$.
The {\em Carnot Caratheodory metric \/}  $d: \HH^n \times \HH^n \to \RR$ is  defined as \be d(p,q)=\inf\{l_H(\Ga)|\Ga \text{ is a horizontal curve from } p \text{ to } q\}.\ee

 A geodesic connecting two points in $\HH^n$ is a horizontal curve of shortest $\HH^n$ length connecting them. It is known that the left translation in $\HH^n$ i.e., $\tau_p : \HH^n \to \HH^n $ defined as $\tau_p(x) = p*x$ is an isometry of $\HH^n$. In particular, the left translation of a geodesic is a geodesic. This helps us visualize the geosdesics of the Heisenberg group better as for all $p,q \in \HH^n$, if $\Ga$ is a geodesic connecting $0$ to $p^{-1}*q$, then $p*\Ga$ is a geodesic connecting $p$ to $q$.  The Heisenberg group is complete, i.e., there exists a geodesic connecting every pair of points in $\HH^n$ (\cite{Zimmerman}). In view of (\ref{segment}), the  straight line segments in the plane $\{t = 0\} \subset \HH^1$ are geodesics  iff they pass through the origin.   For any points $p= (x_1,y_1, 0)$ and $q= (x_2, y_2, 0)$ satisfying $x_2 y_1 - x_1 y_2 \neq 0$, applying translation by $p^{-1}$ to $\HH^1$ so that the new points are $p^{-1} * p = (0,0,0)$ and $ p^{-1} * q = ( x_2-x_1, y_2 - y_1, 2(y_2 x_1 -  x_2 y_1 ) )$  we see that the geodesic joining the translated points does not lie on the plane $\{t =0\}$. Translating this geodesic back by $p$, we will obtain the required geodesic between $p$ and $q$ which is not a straight line segment. For $\HH^1$, a horizontal curve joining the origin and a point $(x, y, t)$ { with $t\neq0$} is a geodesic if and only if it is the lift of circular arc joining the origin with $(x,y,0)$ in the plane $t =0$ whose convex hull has the area equal to scalar multiple of $t$. While, a horizontal curve $\Ga$ from $(0,0,0)$ to $(0,0,T)$ in $\HH^1$ is a geodesic iff the projection of the trace of $\Ga$ in $\RR^2 \times \{0 \}$ is a circle {(see \cite{Capogna}, \cite{Zimmerman})}.  To describe the geodesics in  $\HH^n$ we will use the parametrization  given in{  \cite{Capogna}, \cite{Zimmerman}}.  We also recall the following results from \cite{Zimmerman}, which will be useful for proving our Theorems.
\begin{Pro}\label{zimmer}

{\bf (i)} A horizontal curve $\Ga: [0,1] \to \HH^n = \CC^n \times \RR$ of constant speed connecting $(0,0) \in \CC^n \times \RR$ to $(0,\pm T) \in \CC^n \times \RR$,  $ (T>0)$ is a geodesic iff 
\be  \Ga(s) = \bigg( (1-e^{\mp 2 \pi i s })(A+iB), \pm 2 \| A+iB \|^2 (2 \pi s - \sin 2 \pi s)  \bigg) \label{geodeq} \ee 
where $A=(A_1, A_2,...,A_n), B=(B_1, B_2,...,B_n) \in \RR^n$ are such that
\be \label{AB}
 \| A+iB \|^2 = A_1^2+A_2^2+...+A_n^2 +B_1^2+B_2^2+...+B_n^2 = \frac{T}{4 \pi}.\ee
{\bf (ii)} There are infinitely many geodesics connecting $(0,0,0)$ to $(0,0,T)$. All of them can be obtained from one via a rotation about the $t$ axis.
	
{\bf (iii)} For any $ q \in \HH^n$ which is neither on the $t$-axis nor in the subspace $\CC^n \times \{ 0\}$, there is a unique geodesic connecting the origin to $q$. This geodesic is part of a geodesic connecting the origin to a point on the $t$-axis.

\end{Pro}
To simplify notations, we will choose one particular geodesic  connecting $(0,0)$ to $(0,T)$ denoted by $\Ga_T$ as follows. Since $A+iB $ is any point on the sphere in $\CC^n$  centered at the origin with radius $\frac{|T|}{4 \pi}$, we may choose a particular point 
	\[  A = R =   \left( \dfrac{\sqrt{|T|}}{2\sqrt{n\pi}}, \dfrac{\sqrt{|T|}}{2\sqrt{n\pi}}, \ldots, \dfrac{\sqrt{|T|}}{2\sqrt{n\pi}} \right) \in \RR^n \mbox{~~ and~~} B= 0 \in \RR^n\] 
	so that 
\bea\label{GT} \Ga_T(s) &=& ((1-e^{\mp i s})R,\pm 2\|R\|^2(s -\sin s)) =((1-e^{\mp is})R,
\pm \frac{|T|}{2\pi}(s -\sin s)) \nonumber\\
&& ~~~\mbox{with}~~~ R  = \left( \dfrac{\sqrt{|T|}}{2\sqrt{n\pi}}, \frac{\sqrt{|T|}}{2\sqrt{n\pi}}, \ldots, \dfrac{\sqrt{|T|}}{2\sqrt{n\pi}} \right) \in \RR^n,\, s \in [0, 2\pi].  \label{Ga_T}\eea
We will refer to $\Ga_T$ as a generating geodesic for given $T$ and  in the following, any geodesic between two points on the $t$-axis will be obtained by applying suitable rotation of $\Ga_T$ about the $t$-axis.   
Thus, a Heisenberg bubble is the ``surface of revolution'' obtained by the generating curve $\Ga_T$ and  the whole space $\HH^n$ is foliated by the Heisenberg bubbles. 
Any point $p =(z, t) \in \HH^n$  for which $z\neq 0$ lies on a unique Heisenberg bubble and  a geodesic joining the origin and  the point $p$ is the arc of the unique geodesic passing through $p$ obtained by rotation of $\Ga_T$ for  some $T \in \R$, $T \neq 0$.

\section{Geodetically Convex Sets  in $\HH^n$}  \label{sets}
\setcounter{equation}{0}
Using Heisenberg translation,  any two points on straight line parallel to $t-$axis can always be translated to points on the $t-$axes. Hence, without loss of generality we may assume that the given two points are the origin $ \xi_0= o$ and $ \eta_0=( 0, T_0)$ with $T_0 > 0$.   Let $A=\{o , \eta_0 \} \subset \HH^n$ and let  $\GG(A)$ denote the geodetic convex hull of the set $A$.  Define 
\be \GG^1(A)	 = \cup \{ \Ga(\xi_0,\eta_0): \Ga(\xi_0,\eta_0) \mbox { is a geodesic connecting } \xi_0  \mbox{ and } \eta_0  \}\ee and for $m \geq 2$ the iterated sets
	 		 \be\label{gm}\GG^m(A)=\cup \{ \Ga(\xi,\eta): \Ga(\xi,\eta) \mbox { is a geodesic connecting any two points } \xi, \,  \eta \in \GG^{m-1}(A) \}.  \ee
It is easy to verify that \[ \GG(A) = \cup_{m=1}^\infty \GG^m(A). \]
Note that $\GG^1(A)$   is  the geodesic bubble discussed in Section \ref{pre}. Infact as in case for $\HH^1$, 	 for $\HH^n$ too we see that the iterated sets $\GG^m(A)$, $m \geq 1$ 
	 are never geodetically convex, leading to the conclusion of the Theorem \ref{main01}. 
Precisely, we prove the following theorem in this section. 
\begin{theorem} \label{thm1}   $\cup_{m=1}^\infty \GG^m(A) = \HH^n$.
\end{theorem}

 We begin by listing the properties of the sets $\GG^m(A)$ for any $m \in \NN$. 
\begin{Pro} For all  $m \in \NN$, the set  $\GG^m(A) $ defined in (\ref{gm}) is \\
(i) invariant under the map \[{\mathbb J}: (z,t) \mapsto (\bar z, t);\]
(ii) invariant under the map \[ (x,y,t) \mapsto (x,-y,-t);\]
(iii)	 rotationally invariant about the t axis, i.e., if for $\theta = (\theta_1, ..., \theta_n) \in {\mathbb T}^n$ where $ {\mathbb T}^n =[0, 2\pi]^n$, if the rotation about the $t$-axis   $R_\theta: \HH^n \to \HH^n$ is given by 
\be
 R_\theta (z,t) = (e^{i\theta_1} z_1, \ldots, e^{i\theta_n} z_n,t);
\ee

	then \be R_\theta(\GG^m(A))=\GG^m(A);\ee
(iv) invariant under the reflection about the hyperplane $t=\frac{T_0}{2}$ in $\CC^n \times \RR$.
	\end{Pro}
\proof (i) We note that if $\Ga_{T_0}(s) \in \GG^1(A)$ is a geodesic, then  the curve
\[\si(s) = ((1-e^{is})R, 2\|R\|^2(-s +\sin s))\] is a geodesic joining the origin to the point $-T_0$. Since translation in $t$ variable is an isometry, we see that $({\bf 0}, T_0) * \si(s) $  is a geodesic joining origin and $T_0$ and hence belongs to $\GG^1(A)$.  Similar argument works for all $\GG^m(A)$, $m \geq 2$.

(ii) follows from definition of $\Ga_T$.

 (iii) Since $\GG^1(A)$ is Heisenberg bubble, from definition (\ref{GT}), it follows that $\GG^1(A)$ is rotationally invariant. For any two points $ \xi = (z,t)$, $\eta =(w, s) \in \GG^1(A)$, let $\si[\xi, \eta]$ denote a geodesic joining them. Note that  $R_\theta(\xi)$ and $R_\theta (\eta) \in \GG^1(A)$ for any $\theta \in S^{2n}$, and since $R_\theta $ is an isometry of $\HH^n$, it follows that 
\[ \si[R_\theta(\xi), R_\theta(\eta)] = R_\theta \si[\xi, \eta] \] is a geodesic joining $R_\theta(\xi)$ and $R_\theta (\eta)$. By definition, $R_\theta \si[\xi, \eta] \subset \GG^2(A)$ for any $\theta \in [0, 2\pi]^{n}$.  Since any point  of $ \GG^2(A)$ lies on some geodesic joining any two points of $\GG^1(A)$, the conclusion follows. We can now proceed by induction to deduce the result.

\noindent (iv) We first observe that the Heisenberg bubble $\GG^1(A)$ is symmetric with respect to the plane $t = \frac{T_0}{2}$. For $ \xi_0:=\Ga_{T_0} (s_0) = ((1-e^{-is_0})R, 2\|R\|^2(s -\sin s_0)) $ a given point on $\GG^1(A)$, its reflection with respect to the plane  $t = \frac{T_0}{2}$ is $\hat \xi_0 :=  ((1-e^{-is_0})R, T_0 -  2\|R\|^2(s_0 -\sin s_0))  $. Now from the expression (\ref{GT}) it can be easily verified that the point \Bea
\Ga( 2\pi - s_0) &=&  ((1-e^{-i(2\pi -s_0)})R, 2\|R\|^2((2\pi-s_0)  -\sin (2 \pi - s_0)) \\
&= & ((1-e^{i s_0})R, T_0  - 2\|R\|^2(s_0  - \sin  s_0 ))  \in \GG^1(A).  \Eea
From (i)  we see that 
\Bea {\mathbb J} (\Ga_{T_0}(2\pi - s_0))& =& (\overline{(1-e^{i s_0})}R, T_0 - 2\|R\|^2(s_0  - \sin  s_0) ) \\
&=&  ((1-e^{-i s_0})R, T_0 - 2\|R\|^2(s_0  - \sin  s_0) )  = \hat \xi_0 \in \GG^1(A).   \Eea
\qed 

Let $\tau(p) \in \GG^{m+1}(A)$ where $\tau: \CC^n\times \RR \to \CC^n\times \RR$  defined as $\tau(x+iy,t)=(x+i0,t)$, i.e., \[\tau(x_1, ...,x_n,y_1,...,y_n,t)=(x_1, ...,x_n,0,...,0,t)\] denote the projection of the point $\xi = (z,t)$ onto the space $ \{y= 0\}$.  
		\begin{lemma} \label{proj} For any $m\in \NN$, $ \tau \GG^m(A) \subset \GG^{m+1}(A)$. \end{lemma}
	\begin{proof}
	For a point $\xi = (z_0, t_0) \in \GG^1(A)$,   the straight line segment $\ga(s) = ((1-s) z_0 + s \bar{ z_0}, t_0)$, $s \in [0,1]$  joining $\xi $ and the point ${\mathbb J} \xi $ is a geodesic. Hence by definition of the set $\GG^2(A)$,  
	\[\ga(s) \in  \GG^2(A) \quad \mbox{ for all } s \in [0,1].\]
In fact, since the set $\GG^2(A)$ is rotationally invariant, we conclude that
\[ \{ (z,t_0): |z| \leq |z_0| \} \subset \GG^2(A)\] for any $\xi \in \GG^1(A)$ and hence the solid Heisenberg bubble  is a subset of $\GG^2(A)$. In particular, 
\[\tau \GG^1(A) \subset \GG^2(A).\]
	The proof can now be completed using induction on $m$.
	
\end{proof}
 
 In the following, we further let $T_0 =1$ to simplify notations, so that 
 \[ p =({\bf 0}, 0) \mbox{ and }  q=({\bf 0}, 1) \] and  $\Ga_1:[0, 2\pi] \to \HH^n$ defined as
 \be
 \Ga_1(s) = ((1-e^{-is})R_1, 2(s-\sin s)\|R_1\|^2)  ~~~\mbox{where}~~~ R_1=\left(\frac{\|R_1\|}{\sqrt{n}},..., \frac{\|R_1\|}{\sqrt{n}}\right)\ee
  with $\|R_1\|^2 = \frac{1}{4 \pi}$ is the generating  geodesic  connecting $p$ to $q$     as in (\ref{Ga_T}). Consider the sequence of points $\{p_m, q_m\}$, $m \geq 2$ defined iteratively as 
 \be\label{pq2} p_2 =\tau \bigg(\Ga_{1}(\frac{\pi}{2})\bigg),\, q_2=\tau \bigg(\Ga_{1}(\frac{3 \pi}{2})\bigg)  \ee and 
 \be\label{pqm} p_m =\tau \bigg(\sigma_{m-1}(\frac{\pi}{2})\bigg), \,q_m =\tau \bigg(\sigma_{m-1}(\frac{3 \pi}{2})\bigg), \quad m \geq 2 \ee
where $\sigma_{m-1}$ is the generating geodesic joining $p_{m-1}$ and $q_{m-1}$.
For the proof of Theorem \ref{thm1}, we need to keep track of how the $t$- th coordinate of  grows with $m$ and hence in the following Lemma, we derive the explicit expression for $\Ga_m$, $m \geq2$.
\begin{lemma} \label{Gen Ga} The geodesic joining 
 $ p_2 =\tau \bigg(\Ga_{1}(\frac{\pi}{2})\bigg),\, q_2=\tau \bigg(\Ga_{1}(\frac{3 \pi}{2})\bigg)  \in \GG^2(A) $ is
 \be
 	\sigma_2(s) = \bigg(R_1 +(1-e^{-is})R_2, 2\big(\frac{\pi}{2}-1\big)\|R_1\|^2+2(s-\sin s) \|R_2\|^2 -2 \sin s \langle R_1, R_2 \rangle\bigg)
\label{Ga_2}. \ee
 For each $m \geq 3$,  the points 
 \be p_m =\tau \bigg(\sigma_{m-1}(\frac{\pi}{2})\bigg), \,q_m =\tau \bigg(\sigma_{m-1}(\frac{3 \pi}{2})\bigg) \in \GG^{2m-2} (A) \ee 
 and  the geodesic $\sigma_m: [0,2 \pi] \to \HH^n$ connecting  $p_m$ to $q_m$ is  given by 
 \begin{align}
\sigma_m(s) =& \bigg(\sum_{j=1}^{m-1} R_j +(1-e^{-is})R_m, t_m(s) \bigg) \mbox{ where } \nonumber \\
t_m(s) = & 2\big(\frac{\pi}{2}-1\big)\big(\sum_{j=1}^{m-1} \|R_j\|^2\big)+2(s-\sin s) \|R_m\|^2  
-2 \sum\limits_{j=2}^{m-1}( \langle  \sum\limits_{k=1}^{j-1} R_k, R_j \rangle )
-2 \sin s \langle \sum_{j=1}^{m-1} R_j , R_m \label{Ga_m}\rangle. 
\end{align}
Here $\langle \cdot, \cdot \rangle$ denotes the usual (real) inner product and 
\[R_m =\bigg(\frac{\|R_m\|}{\sqrt{n}},...,\frac{\|R_m\|}{\sqrt{n}}\bigg)\] with
\[
\|R_m\|^2 = \frac{1}{4 \pi}\big[ 2 \| R_{m-1}\|^2 (\pi+2) +2\langle R_1 +...+R_{m-2}, R_{m-1}\rangle \big]. \]

\end{lemma}
\begin{proof}

From Lemma \ref{proj}, we observe that 
	\begin{align}
	p_2 &=\tau \bigg(\Ga_{1}(\frac{\pi}{2})\bigg)=\big(R_1, 2 (\frac{\pi}{2}-1)\|R_1\|^2\big)\in \GG^2(A),\\
	q_2 &=\tau \bigg(\Ga_{1}(\frac{3\pi}{2})\bigg)=\big(R_1, 2 (\frac{3\pi}{2}+1)\|R_1\|^2 \big) \in \GG^2(A) \\ \mbox{with } R_1&=\left(\frac{\|R_1\|}{\sqrt{n}},..., \frac{\|R_1\|}{\sqrt{n}}\right) \mbox{ and }  \|R_1\|^2 = \frac{1}{4 \pi}. 
	\end{align}
Both $p_2$ and $q_2$ lie on a axis parallel to the $t-$axis.  Hence, a  geodesic $\Ga_2: [0, 2 \pi] \to \HH^n$ in $\GG^3(A)$ connecting $p_2$ to $q_2$ can be obtained by translating the generating geodesic joining the origin and $-p_2*q_2$ by $p_2$, i.e.,
	\[ \sigma_2(s) = p_2 * \Ga_2(s)\]
	where $\Ga_2:[0,2 \pi] \to \HH^n$, as in (\ref{GT}) is a geodesic connecting $(0,0)$ to $-p_2*q_2$ with \[R_2 = \left(\frac{\|R_2\|}{\sqrt{n}},..., \frac{\|R_2\|}{\sqrt{n}}\right) \mbox{  and  } \|R_2\|^2 = \frac{1}{4 \pi}\big[2\|R_1\|^2(\pi + 2)\big].\] Therefore, we have 
	\begin{align}
	\Ga_2(s) =& \bigg(R_1 +(1-e^{-is})R_2, 2\big(\frac{\pi}{2}-1\big)\|R_1\|^2+2(s-\sin s) \|R_2\|^2 -2 \sin s \langle R_1, R_2 \rangle\bigg)
	\end{align} which is as in (\ref{Ga_2}).
	
		We prove (\ref{Ga_m}) by induction on $m\geq3$.	For $m=3$,  again using the Lemma \ref{proj}, we observe that 
			\begin{align}
		p_3 &=\tau \bigg(\sigma_{2}(\frac{\pi}{2})\bigg)=\bigg(R_1 +R_2, 2\big(\frac{\pi}{2}-1\big)\big(\sum_{j=1}^{2} \|R_j\|^2\big) -2 \langle R_1, R_2 \rangle\bigg) \in \GG^4(A),\\
		q_3 &=\tau \bigg(\sigma_{2}(\frac{3\pi}{2})\bigg)\nonumber \\&=\big(R_1+R_2, 2 (\frac{\pi}{2}-1)\|R_1\|^2 +2 (\frac{3\pi}{2}+1)\|R_2\|^2 +2 \langle R_1, R_2 \rangle \big) \in \GG^4(A) . 
		\end{align}
		A geodesic $\sigma_3: [0, 2 \pi] \to \HH^n$ in $\GG^5(A)$ connecting $p_3$ to $q_3$ can be given by 
\be \sigma_3(s) = p_3 * \Ga_3(s)\ee
where $\Ga_3:[0,2 \pi] \to \HH^n$, as in (\ref{GT}) is a geodesic connecting $(0,0)$ to $-p_3*q_3$ with \be R_3 = (\frac{\|R_3\|}{\sqrt{n}},..., \frac{\|R_3\|}{\sqrt{n}}) \mbox{  and } \|R_3\|^2 = \frac{1}{4 \pi}\big[2\|R_2\|^2(\pi + 2) +4 \langle R_1, R_2 \rangle\big]. \ee
 Therefore, we have 
\begin{align}
\sigma_3(s) =& \bigg(\sum_{j=1}^{2} R_j +(1-e^{-is})R_3, 2\big(\frac{\pi}{2}-1\big)\big(\sum_{j=1}^{2} \|R_j\|^2\big)+2(s-\sin s) \|R_3\|^2 \nonumber \\
&-2 \langle R_1, R_2 \rangle -2 \sin s \langle R_1+R_{2}, R_3 \label{Ga_3}\rangle\bigg) 
\end{align} which is as in (\ref{Ga_m}), proving the result for $m=3$.

By induction, assuming  that the expresssion for $\sigma_m$, $m\geq 3$ is given by (\ref{Ga_m})  we  derive it for  $\sigma_{m+1}$.
	Since, $p_m, q_m \in \GG^{2m-2}$, therefore, by the definition of $\GG^n(A)$ we can say, $\Ga_m$ lies in $\GG^{2m-1}(A)$. Now using Lemma \ref{proj} again, we see that 
		\begin{align}
p_{m+1} =&\tau \bigg(\Ga_{m}(\frac{\pi}{2})\bigg)\nonumber \\
=&\bigg(\sum_{j=1}^{m} R_j , 2\big(\frac{\pi}{2}-1\big)\big(\sum_{j=1}^{m} \|R_j\|^2\big)-2 \langle R_1, R_2 \rangle -2\langle R_1 + R_2, R_3 \rangle-\nonumber \\
&\;\;\;\;...-2\langle R_1+...+R_{m-2}, R_{m-1} \rangle -2 \langle R_1+...+R_{m-1}, R_m \rangle\bigg) \in \GG^{2m}(A),\\
q_{m+1} =&\tau \bigg(\Ga_{m}(\frac{3\pi}{2})\bigg) \nonumber \\
=&\bigg(\sum_{j=1}^{m} R_j , 2\big(\frac{\pi}{2}-1\big)\big(\sum_{j=1}^{m-1} \|R_j\|^2\big)+2\big(\frac{3\pi}{2}+1\big)\|R_m\|^2-2 \langle R_1, R_2 \rangle -2\langle R_1 + R_2, R_3 \rangle -\nonumber \\ &
...-2\langle R_1+...+R_{m-2}, R_{m-1} \rangle +2 \langle R_1+...+R_{m-1}, R_m \rangle\bigg) \in \GG^{2m}(A).
\end{align}
So, \be -p_{m+1}*q_{m+1} = (0, 2 \| R_{m}\|^2 (\pi+2) +4\langle R_1 +...+R_{m-1}, R_{m}\rangle).\ee

	A geodesic $\sigma_{m+1}: [0, 2 \pi] \to \HH^n$ connecting $p_{m+1}$ to $q_{m+1}$ can be given by 
\be \sigma_{m+1}(s) = p_{m+1} * \Ga_{m+1}(s)\ee
where $\Ga_{m+1}:[0,2 \pi] \to \HH^n$, as in (\ref{GT}) a geodesic connecting $(0,0)$ to $-p_{m+1}*q_{m+1}$ with 
\be R_{m+1} = \left(\frac{\|R_{m+1}\|}{\sqrt{n}},..., \frac{\|R_{m+1}\|}{\sqrt{n}}\right)\ee and \be \|R_{m+1}\|^2 = \frac{1}{4 \pi}\big[2 \| R_{m}\|^2 (\pi+2) +4\langle R_1 +...+R_{m-1}, R_{m}\rangle \big]. \nonumber \ee Therefore, we have 
\begin{align}
\sigma_{m+1}(s) =& \bigg(\sum_{j=1}^{m} R_j +(1-e^{-is})R_{m+1}, 2\big(\frac{\pi}{2}-1\big)\big(\sum_{j=1}^{m} \|R_j\|^2\big)+2(s-\sin s) \|R_{m+1}\|^2 \nonumber \\
&-2 \langle R_1, R_2 \rangle -2\langle R_1 + R_2, R_3 \rangle-...-2\langle R_1+...+R_{m-1}, R_{m} \rangle \nonumber  \\
&-2 \sin s \langle R_1+...+R_{m}, R_{m+1} \rangle\bigg)
\end{align} which is of the form (\ref{Ga_m}).
\end{proof}

{\bf Proof of Theorem \ref{thm1}}
Since the height of geodesic bubble determines its size, to prove the Theorem \ref{thm1} it suffices to show that 
\be || R_m|| \to \infty \quad \mbox{as} \quad m \to \infty .\ee
From Lemma \ref{Gen Ga}, we have for all  $m>3$,
		\begin{align}
	\|R_m\|^2&= \frac{1}{4 \pi}\big[2 \|R_{m-1}\|^2(\pi+2)+4 \langle R_1+...+R_{m-2}, R_{m-1} \rangle\big] \nonumber \\
	&= \|R_{m-1}\|^2\big(\frac{\pi+2}{2 \pi} \big)+ \frac{1}{\pi}\langle R_1+...+R_{m-2}, R_{m-1} \rangle \nonumber
	\end{align}
	where each \be R_k = \big( \frac{\|R_k \|}{\sqrt{n}}, ...,  \frac{\|R_k \|}{\sqrt{n}}\big) \text{ for } k = 1,2,..., m. \nonumber\ee
Combining the above equations gives us for all  $m >3 $, 
\be\label{00} \|R_m\|^2 = \|R_{m-1}\|^2\big(\frac{\pi+2}{2 \pi} \big)+ \frac{1}{\pi}( \| R_1\| +...+\| R_{m-2}\|) \| R_{m-1} \| .\ee
	We claim that  \be \label{main1} \|R_m\|^2 > \|R_{m-1}\|^2 \big(\frac{\pi+2}{2 \pi} \big) + \frac{1}{5}\|R_{m-1}\|^2 \;\; \mbox{ for all } m >3. \ee 

To prove the claim, in view of  (\ref{00}) it is enough to prove
 \[ 5 \big(\| R_1\|+\|R_2\|+...+\|R_{m-2} \| \big) > \pi \| R_{m-1} \| \; \forall m > 3. \]
Since, $\| R_j \| $s are positive, this is equivalent to showing
\be\label{01} 25 \big(\| R_1\|+\| R_2\|+...+\|R_{m-2} \| \big)^2 > \pi^2 \| R_{m-1} \|^2  \;\;\; \mbox{ for all } m > 3,  \ee
Substituting the value of $||R_{m-1}||^2$, we have RHS of (\ref{01}) is
\[  \pi \| R_{m-2} \| ^2 \frac{(\pi+2)}{2}+ \pi \big(\| R_1 \|+ \| R_2 \|+...+ \|R_{m-3} \| \big) \| R_{m-2} \|.   \]
While 
	by computation, we see that the LHS of (\ref{01}) is 
\Bea &&25 \big( \| R_1\|+\|R_2\|+...+\|R_{m-2}\| \big)^2\\
&& = 25 (\| R_{m-2}\|^2 + 2\| R_1 \| \| R_{m-2} \| +2\| R_2\| \| R_{m-2}\| +...+2\| R_{m-3}\| \| R_{m-2}\| +\text{ other positive terms} ).\Eea
Since, 
\be 25 \| R_{m-2} \| ^2> \pi \| R_{m-2}\| ^2 \frac{(\pi+2)}{2}, \nonumber \ee
\be 50 \| R_1\| \| R_{m-2}\| > \pi \| R_1\| \| R_{m-2} \|,  \nonumber \ee
\be :  \nonumber \ee
\be 50 \| R_{m-3}\| R_{m-2} \| >  \pi \| R_{m-3}\| \| R_{m-2} \|, \nonumber \ee
(\ref{01})  and hence (\ref{main1}) is true. This gives us \be \label{main4} \| R_m\| ^2 > (1+c)\| R_{m-1}\| ^2 \;\; \forall m >2   \ee
where $c=\frac{10-3\pi}{10 \pi} > 0$. It follows that $\lim\limits_{m\to \infty} \|R_m \| ^2 = \infty$ for otherwise if $\lim\limits_{m\to \infty} \| R_m \|^2 = \alpha< \infty$ then taking limit as $m \to \infty$ in (\ref{main4}) we get a contradiction. 

\qed

For the proof of Theorem \ref{thm2}, we require  the extensions Lemma 3.1 in \cite{Monti} to higher dimensions. Though the steps are similar, we have included the proof for the sake of completeness. 

\begin{lemma}\label{graph}
Let $\Om \subset \HH^n$. Suppose there exists a point $\xi_0 \in \Om$  and a neighbourhood of $\xi_0$   in $\Om$ which can be expressed as a graph of a continuous function over the  $z$-plane. Then $\Om$ is not geodetically convex.
\end{lemma}
\begin{proof} 
Without loss of generality, let $\xi_0 = o$, the origin and suppose by contradiction that $\Om$ is geodetically convex. If $\Om$ contains two distinct points on a line parallel to $t$ axis then by Theorem \ref{main01},  $\Om = \HH^n$ and there is nothing to prove. If not, then any two points in $\xi,\, \eta \in \Om$ can be joined by a unique geodesic denoted by $\ga_{\xi\eta}$, say.
Let $U \subset \R^{2n}$ be a neighbourhood of origin and $f: U \to \RR$  be a  continuous function with $f(0) =  0 \in \R$ such that \be V :=
 \{ (x, y, f(x,y)) = (z, f(z)) : (x,y) \in U \} \subset \Om.\ee
 Choose $r_0 > 0$ sufficiently small such that $B(0, r_0) \subset \subset U$ and define \[ g(\omega) = f( \frac{\om}{\sqrt n}, \frac{\om}{\sqrt n}, \ldots,  \frac{\om}{\sqrt n} ) - f (-  \frac{\om}{\sqrt n}, -  \frac{\om}{\sqrt n}, \ldots, - \frac{\om }{\sqrt n}) \mbox{ for } \om \in S(r_0) \subset \CC\] 
where $S(r_0)$ is a circle in the complex plane centered at the origin with radius $r_0$. We claim that there exists $\om_0  \in  S(r_0)$ such that 
 \be g(\om_0) = 0. \ee
 To see this, we note that $g$ is a non trivial continuous odd function defined on a circle, and hence there exists $\om_0\in S( r_0)$ such that $g(\om_0) = 0$, i.e.,
 \be
 f( \frac{\om_0}{\sqrt n}, \frac{\om_0}{\sqrt n}, \ldots,  \frac{\om_0}{\sqrt n} ) = f (-  \frac{\om_0}{\sqrt n}, -  \frac{\om_0}{\sqrt n}, \ldots, - \frac{\om_0 }{\sqrt n}) = t_0 \mbox{ (say)}. 
 \ee
 Note that 
 \[ z^0 = ( \frac{\om_0}{\sqrt n}, \frac{\om_0}{\sqrt n}, \ldots,  \frac{\om_0}{\sqrt n} ) \in \pa B(0, r_0) \subset U. \] As $\Om$ is geodetically convex, the geodesic (line segment)  $ \ga(s) =(  sz^0, t_0 ) \in \Om$ and infact $ sz^0 \in U$ for all $s \in [-1, 1]$.  
 
 If $t_0 \neq 0$ then already $(0, t_0)$ and origin are in $\Om$ and $\Om = \HH^n$. Hence we must have $t_0 = 0$ and $(s z^0, 0) \in V$ for all $ s \in [-1,1]$. Let $L := \{ (sz^0,0): s \in [0,1]\}$ denote the line segment in $V \subset \Om$.
 
 Now we will show that there exists two points $ p$ and $ q \in B(0, r_0)$ such that the geodesic joining these two points $\tilde\xi = (p, f(p))$ and $\tilde\eta = (q, f(q)) \in \Om$ projects onto an arc intersecting the line segment $L$, again a contradiction. Without loss of generality, we may assume
 \be \left. \begin{array}{lll}  z^0  &= &( x_1^0,\ldots,x_n^0, 0, \ldots, 0); \\
  \tilde \xi & =& (p, \ldots, p,\var,\ldots ,\var, f(p, \ldots, p,\var,\ldots ,\var) ) =  (p, \ldots, p,\var,\ldots ,\var, t_1 )\\
 \tilde \eta &= & (p, \ldots, p, -\var,\ldots ,-\var, f(p, \ldots, p,-\var,\ldots ,-\var)) = (p, \ldots, p,-\var,\ldots ,-\var, t_0 ),
    \end{array}\right\}\ee
where $0< p, \var \in \RR$ such that  $(p, \ldots, p, 0)$ lies sufficiently close to the origin on the line segment $L$ and $\var >0$ is small to be chosen later. It can be verified that the geodesic joining $\tilde\xi$ and $\tilde\eta$ is given by the curve
 \[\sigma(s) = (x(s), y(s), t(s)) , \quad s \in [0, s_0] \subset [0, 2\pi) \] where \bea x_i(s) & = & p+  \var \left(   \sin s_0 \frac{(1-\cos s)}{(1-\cos s_0)} -  \sin s  \right), \, 1 \leq i \leq n,\\
 y_i(s) & = & -\var +  \var \left( (1- \cos s) + \frac{\sin s_0 \sin s}{(1-\cos s_0)} \right), \,  1 \leq i \leq n, \\
 t(s) & =  &    t_0 + \frac{ 4n \var^2}{(1-\cos s_0)} (s- \sin s) + 2 \sum_{i = 1}^n (-\var x_i(s) - p y_i(s))  \\
 && s \in [0, s_0]  \mbox{ such that  } \frac{\var(s_0 -\sin s_0)}{(1-\cos s_0)} = p,  \, \frac{ 4n \var^2}{(1-\cos s_0)} = t_1 -t_0 + 2p \var.
 \eea
 Observe that $t(s) \neq 0$ for all $s \in (0, s_0)$ and hence the geodesic $\sigma$  does not lie in the hyperplane $t =0$. For $p$ close to origin and $\var$ sufficiently small, $\si(s)$  projects onto the curve  $(x(s), y(s))$ 
which intersects the line segment $L$. Thus we have found two distinct points in $\Om$ which are parallel to the $t$ axis which is a contradiction to $V$ being a graph.

\end{proof}

\subsection{Proof of Theorem \ref{thm00}}
Let $o$, $\xi$ and $\eta$ be three distinct points in a set $A$, all of do not lie on a geodesic. Without loss of generality, assume that $o$ is origin. We prove that the geodetic convex hull ${\mathcal G}(A) = \HH^n$. 

Due to Theorem \ref{main01}, if $A $ contains two distinct points which lie on a line parallel to the $t$-axis then ${\mathcal G}(A) = \HH^n$ and there is nothing to prove. Hence suppose that such points are not in $A$ and that every pair of  points in ${\mathcal  G}(A)$ can be joined by a unique geodesic. 

Let $d(\xi, \eta) = d_{\xi\eta}$ denote the distance between the points $\xi$ and $\eta$
 and $\sigma  : [0, d_{\xi\eta}] \to \HH^n$ denote the unique geodesics joining them with $\si(0) = \xi$ and $\si(d_{\xi\eta}) = \eta$. Also let $d(o, \si(s)) := d_{\si(s)}$ denote the distance between the origin $o$ and  $\si(s)$ as $s$ varies in $[0,  d_{\xi\eta}]$ and define  
$\ga_{\si(s)}: [0, d_{\si(s)}] \to \HH^n$ to be the (unique) geodesic joining origin to $\si(s)$. Arguing as in 2. of proof of Theorem 1.2 (pg.194, \cite{Monti}), we see that all of them are distinct and do not intersect each other.   
Consider the open set  $U := \{(\theta, s) : \theta \in (0,d_{\si(s)}), \, s \in (0, d_{\xi\eta})\}\subset \RR^{2}$ and define the map
 $\Phi: U \to \HH^n $ as
\[ \Phi( \theta, s) = \ga_{\si(s)}(\theta) . \]
Then $\Phi$ is continuous. Let $ \Phi\left( U \right) =: V$. Then  every $\nu \in V$ lies on a unique geodesic $\ga_{\si(s)} $ such that  $ \ga_{\si(s)}(\theta) = \nu$ 
and we define $P: V \to \RR^{2n}$ as 
\be P (\nu) = \pi_{xy}\ga_{\si(s)}(\theta) \in \RR^{2n},\ee
 where $\pi_{xy}$ denotes  the projection onto the $ x-y $ i.e., $z$- hyperplane of the unique geodesic  $\ga_{\si(s)} $.
The map $P$ is continuous and one-one. For  if $\nu_1$, $\nu_2 \in V$ such that 
$P(\nu_1) = P(\nu_2)$, then $\pi_{xy}\ga_{\si(s_1)}(\theta_1) = \pi_{xy}\ga_{\si(s_2)}(\theta_2) = (x_1^0,\ldots, x_n^0, y_1^0, \ldots, y_n^0) \in \RR^{2n} $ say. 
But then $\nu_1$ and $\nu_2$ lie on a line parallel to $t$- axis and  we conclude that  either  $A = \HH^n$ due to Theorem \ref{main01} or $\nu_1 = \nu_2$. 

Suppose $A \neq \HH^n$, so that $P$ is a one one continuous map. It follows from the Theorem 4.17 of \cite{rudin} that $P$ is a homeomorphism from $V$ onto $ P\left( V \right) =: \tilde V \subset \RR^{2n}$. Define $f:\tilde V \to \RR$ as 
\[ f(x,y) := \pi_t(P^{-1}(x, y)) \] where $\pi_t$ is the projection onto the $t$-th coordinate. Then $V$ is expressed as a graph $(x, y, f(x, y))$ for $(x,y) \in \tilde V$ and we again get contradiction due to Lemma  \ref{graph}. This completes the proof of Theorem \ref{thm00}.

\qed


\section{Geodetically Convex Functions}\label{functions}
Recall that   a function $f:[a,b] \to \RR$ is said to be \textbf{convex} if $\forall x_1, x_2 \in [a,b]$ and $\forall \la \in (0,1)$, we have
	\[ f ( \la x_1 +(1-\la) x_2) \leq \la  f( x_1) +(1-\la) f(x_2). \] In particular,  if $f$ is convex then $\forall x \in [a,b]$,
	\[f(x) \leq f(a) \text{ or } f(x) \leq f(b)\]
since  $\forall x \in [a,b],$
	\[ f(x) \leq m(x-a)+ f(a)\] where $m=\big(\frac{f(b)-f(a)}{b-a}\big)$ is the slope of the line segment joining $(a,f(a))$ and $(b,f(b))$.
	If $m \geq 0$ then  $f(x) \leq f(b)$. Similarly, if $m \leq 0$ then $f(x) \leq f(a)$.

\begin{lemma}
	Let $u: \HH^n \to \RR$ be geodetically convex and $p_0 \in \HH^n$ then the set \be S=\{ p \in \HH^n | u(p)< u(p_0)\} \nonumber \ee is geodetically convex. 
\end{lemma}
\begin{proof}
	Consider any $p_1, p_2 \in S$ and any geodesic $\ga: [0,1] \to \HH^n$ such that $\ga(0)=p_1$ and  $\ga(1)=p_2$.	We claim that  
	\be \ga(s) \in S \;\; \forall s\in [0,1].\ee
	This can be achieved if we show that $u(\ga(s)) < u(p_0) \;\; \forall s \in [0,1]$.
	Now, \be p_1, p_2 \in S \implies u(p_1), u(p_2) < u(p_0). \ee
	Since, $u \circ \ga $ is convex, $\forall s \in [0,1]$, we have \be u \circ \ga (s) \leq u \circ \ga (0) \text{ or }  u \circ \ga (s) \leq u \circ \ga (1).\ee
	In either case, we will have $u \circ \ga (s) < u(p_0) \;\; \forall s \in [0,1]$.
	Hence, proved.
	
\end{proof}
\textbf{Proof of Theorem \ref{thm3.2}:} It suffices to prove that \be u(p)=u(0)\;\; \forall p \in \HH^n. \label{3.7}\ee
Suppose (\ref{3.7}) is not true. Then there exists some $q_0 \in \HH^n$ such that $u(q_0) \neq u(0)$. 
Without loss of generality, let us assume $u(q_0) < u(0)$. 
Define \be S_0 = \{ q \in \HH^n | u(q)< u(0)\}. \ee
Then $S_0 \neq \emptyset$ open subset of $\HH^n$ as $q_0 \in S$ and $u$ is continuous. Therefore, there exists $r>0$ such that $B(q_0, r) \subset  S_0$ and for sufficiently small $\var >0$ we have 
 \be q_0-\epsilon e_{n+1}, q_0+\epsilon e_{n+1} \in B(q_0,r)   \ee where $e_{n+1}=(0,...,0,1)\in \HH^n$. 
By the previous Lemma, $S_0$ is a geodetically convex set containing two points $q_0-\epsilon e_{n+1}, q_0+\epsilon e_{n+1} \in \CC^n \times \RR$ with the same projection in $\CC^n \times \{ 0 \}$. Now using Theorem \ref{thm1}, $S_0=\HH^n$ which is a contradiction as $0\notin S_0$ by definition of $S_0$.

Note that if we assume $u(q_0) > u(0)$ then we define \be S_{q_0} = \{ q \in \HH^n | u(q)< u(q_0)\} \nonumber \ee which is also non empty as $0 \in S_{q_0}$ and similar arguments can be made to arrive at a contradiction.
Hence, our supposition is not true, implying $u$ is a constant function. 

\qed

\end{document}